\magnification 1200
\baselineskip 14pt
\parskip=3pt plus1pt minus.5pt
 2
\font\large=cmbx10 scaled \magstep 2
 4

 1
\input amssym.def

\def\qed{\hfill\vbox{\hrule\hbox{\vrule\kern3pt
\vbox{\kern6pt}\kern3pt\vrule}\hrule}\bigskip}
\def\mapright#1{\smash{
   \mathop{\longrightarrow}\limits^{#1}}}

\centerline{\large A home-made 
Hartshorne-Serre correspondence}
\centerline{Enrique Arrondo
\footnote{(*)}{Supported in part by the funds
of the Spanish research project BFM2003-03971}}
\bigskip
\bigskip

\item{}{\bf Abstract:} We provide an
elementary proof of the Hartshorne-Serre
correspondence for constructing vector bundles
from local complete intersection subschemes of
codimension two. This will be done, as in the
correspondence of hypersurfaces and line
bundles, by patching together local
determinantal equations in order to produce
sections of a vector bundle.

\bigskip
\noindent{\bf Introduction.}
\medskip

It is well-known that a hypersurface of a
smooth algebraic variety can be obtained (in a
unique way) as the zero locus of a section of a
line bundle. In fact the construction
of the line bundle and its section can be done
in a very elementary way, by patching local
equations and it can be taught in any first
course of algebraic geometry. 

If instead one considers subvarieties of
codimension bigger than one, the situation is
very different and only well understood in
codimension two (for some results in
codimension three, see [13], [4] and [5]).
More precisely, in [7], Hartshorne, inspired by
previous works of Serre and Horrocks ([11] and
[9]), proved that a codimension two subvariety
of ${\Bbb P}^n$ is the zero locus of a rank-two
vector bundle over ${\Bbb P}^n$ if and only if
the subvariety is subcanonical (which can be
interpreted as saying that the determinant of
its normal bundle extends to a line bundle $L$
on ${\Bbb P}^n$). This result was independently
proved by Barth and Van de Ven ([3]), and
generalized by Grauert and M\"ulich ([6]) to any
ambient space (in which case the vanishing of
the second order cohomology of $L^*$ is needed).
Finally, following the original technique of
Hartshorne, Vogelaar ([12]) gave the most general
result, proving that any local complete
intersection subscheme of codimension two of a
smooth variety $X$ can be obtained as the
dependency locus of $r-1$ sections of a rank $r$
vector bundle over $X$ of determinant $L$ if and
only if the determinant of its normal bundle
twisted with $L^*$ is generated by $r-1$ global
sections (provided again the vanishing of the
second order cohomology of $L^*$). In both [6]
and [12], the uniqueness of the vector bundle is
obtained provided the vanishing of the first
order cohomology of $L^*$.

However, although this construction (already
known as Hartshorne-Serre correspondence) is
very well-known and thoroughly used, it is very
difficult to provide a good reference of it.
Indeed the general result is only in Vogelaar's
PhD thesis, which is not published
elsewhere, and hence it is usually embarrassing
to use as a reference. Even in the subcanonical
case, although the technique of [6] works in
general, it is written only for the particular
case of projective spaces (as it happens for
all the other proofs I know of this case).

On the other hand, it is also very annoying
that, while the case of codimension one is so
easy to explain to even an undergraduate
student, the techniques for the case of
codimension two are too sophisticated, using in
an essential way the spectral sequence of
local and global Ext). Only in [6] (which is
written in german) there is a more elementary
proof of the subcanonical case. 

The goal of this paper is hence double. On one
hand, we want to provide a reference for the
general Hartshorne-Serre correspondence. On the
other hand, in order to present some new
material, we will give a quite elementary proof
of the main result, namely patching together
local representations of the sections of the
vector bundle we are looking for (hence
imitating the standard proof for codimension
one). This is in fact the method used in [6],
without much details, in the subcanonical case (I
thank the lovely kindness of Sof\'{\i}a Cobo, who
translated for me that paper, so that I learned
that my first draft [2] contained essentially
the proof of [6]). Anyway, the general case
still requires some new tricky ideas that we
develop in this paper. We also hope that our
approach could be extended to other contexts
different from algebraic geometry, and maybe
give also some idea about how to extend this
kind of results to higher codimension. 

In a first section, we will recall the main
result (Theorem 1), and for the sake of
completeness we will also recall its standard
algebraic proof. This will be the only part in
which a good background of algebraic geometry (at
the level of [8]) will be required. For the
rest of the paper, we hope that it will be
readable for a wide range of mathematicians (it
will not be important at all to know what a
scheme is!!!). In the second section, we will
present the minimal background needed to follow
the paper.

The proof of Theorem 1 will be divided in the
next three sections. In section 3, we will
study the main properties that we will require
to an open covering of our general ambient
variety. In section 4, we will discuss how to
construct the $r-1$ sections of the vector
bundle we are looking for. We will eventually
finish the proof of Theorem 1 in section 5, in
which we will see how the cohomological
conditions on $L^*$ imply the existence and
uniqueness of the vector bundle. Finally, we
include a last section for some remarks on
possible generalizations of the result.

Although I usually do not like to do it, I had
to sacrifice the ``pedagogical'' presentation of
the material by the sake of the rigor. In other
words, I decided to avoid sentences like ``we
could have refined our covering so that...'' or
``changing our definition of [...] we can
assume...''. As a result, several definitions
and notations that a priori seem artificial can
only be understood a posteriori. For example,
the strange sign in Lemma 2 is explained after
Lemma 4 (see Remark 5), and the apparently
complicated way of writing the matrices in
section 4 makes sense only in section 5.

\bigskip
\noindent{\bf 1. Statement and the standard
approach.} 
\medskip

Let $X$ be a smooth algebraic variety over an
algebraically closed field $k$. Let $Y$ be a
codimension two subscheme of $X$. We will
denote by ${\cal J}$ the ideal sheaf of $Y$ in
$X$. If we assume that $Y$ is a local complete
intersection, then the conormal sheaf
$N^*:={\cal J}\otimes{\cal O}_Y$ is locally
free, so that we will regard its dual $N$ as a
vector bundle. Assume that $Y$ is the dependency
locus of $r-1$ sections
$\alpha_1,\ldots,\alpha_{r-1}$ of a rank
$r$ vector bundle $E$ over $X$ with $\bigwedge^rE=L$.
This produces an exact sequence 
$$0\to(r-1){\cal
O}_X\mapright{(\alpha_1,\ldots,\alpha_{r-1})}E\to
{\cal J}\otimes L\to 0. \eqno{(1)}$$  Its
restriction to $Y$ produces a long exact
sequence 
$$0\to\bigwedge^2N^*\otimes L_{|Y}\to(r - 1){\cal
O}_Y\mapright{({\alpha_1}_{|Y},\ldots,{\alpha_{r-1}}_{|Y})}
E_{|Y}\to N^*\otimes L_{|Y}\to 0
\eqno{(2)}$$  in which we find out that the kernel
of the middle map is $\bigwedge^2N^*\otimes L_{|Y}$
by just looking at the first Chern classes in
the sequence. Dualizing the first map in (2) we 
obtain in particular that the line bundle
$\bigwedge^2N\otimes L^*_{|Y}$ is generated by $r-1$
global sections $s_1,\ldots,s_{r-1}$ that also
satisfy 
$$s_1{\alpha_1}_{|Y}+\ldots+s_{r-1}{\alpha_{r-1}}_{|Y}=0.$$ 
Hartshorne-Serre correspondence consists of
reversing this process. More precisely:

\proclaim Theorem 1. Let $X$ be a smooth algebraic
variety and let $Y$ be a local complete
intersection subscheme of codimension two in $X$.
Let $N$ be the normal bundle of $Y$ in $X$ and let
$L$ be a line bundle on $X$ such that
$H^2(X,L^*)=0$. Assume that $\bigwedge^2N\otimes
L^*_{|Y}$ has $r-1$ generating global sections
$s_1,\ldots,s_{r-1}$. Then there exists a
rank $r$ vector bundle $E$ over X such that: 
\item{(i)} $\bigwedge^r E=L$; 
\item{(ii)} $E$ has $r-1$ global sections
$\alpha_1,\ldots,\alpha_{r-1}$ whose dependency
locus is $Y$ and such that
$s_1{\alpha_1}_{|Y}+\ldots+s_{r-1}{\alpha_{r-1}}_{|Y}=0$. 
\smallskip\noindent Moreover, if
$H^1(X,L^*)= 0$, conditions (i) and (ii) determine
$E$ up to isomorphism.

The main idea for the standard algebraic proof
is to obtain $E$ as an extension like (1), i.e.
as a suitable element in ${\rm Ext}^1({\cal
J}\otimes L,(r-1){\cal O}_X)$. For this, one
first considers the spectral sequence 
$$E^{p,q}_2:=H^p({\cal E}xt^q({\cal J}\otimes
L,(r-1){\cal O}_X))\Rightarrow E_n:={\rm
Ext}^n(L\otimes{\cal J},(r-1){\cal O}_X))$$
(see [1] Proposition (2.4)).
Then the exact sequence
$$0\to E_2^{1,0}\to E^1\to E_2^{0,1}\to
E_2^{2,0}$$
(see for instance [10] Theorem 11.43), becomes,
under natural identifications:
$$\matrix{0\to H^1(X,(r-1)L^*)\to {\rm
Ext}^1({\cal J}\otimes L,(r-1){\cal
O}_X))\mapright{\varphi}\cr\cr
\mapright{\varphi}{\rm
Hom}((r-1){\cal O}_X,\bigwedge^2N\otimes
L_{|Y}^*)\mapright{\psi}
H^2(X,(r-1)L^*)}\eqno{(3)}$$
Under the hypothesis $H^2(X,L^*)=0$, the map
$\varphi$ is surjective, and hence the element
$\eta\in{\rm Hom}((r-1){\cal
O}_X,\bigwedge^2N\otimes L_{|Y}^*)$
corresponding to the choice of
$s_1,\ldots,s_{r-1}$ will produce an extension as
in (1) (which will be unique if
$H^1(X,L^*)=0$). Hence it is enough to check
that $E$ is a locally free sheaf. The proof
given in [12] is not clear to us, so that we
outline here another (standard) one. We need to
show that ${\cal E}xt^i(E,{\cal O}_X)=0$ for all
$i>0$, and this can be done by applying the
functor ${\cal H}om(\underline{\hskip
0.2cm},{\cal O}_X)$ to the exact sequence (1) we
just constructed. The only difficulty is to show
the vanishing of
${\cal E}xt^1(E,{\cal O}_X)$, but this follows
from the fact that in the exact sequence
$${\cal H}om((r-1){\cal O}_X,{\cal O}_X)\to
{\cal E}xt^1({\cal J}\otimes L,{\cal
O}_X)\to {\cal E}xt^1(E,{\cal O}_X)\to 0$$
the first morphism is canonically identified
with the surjection $\eta:(r-1){\cal
O}_X\to\bigwedge^2N\otimes L_{|Y}^*$ induced by
$s_1,\ldots,s_{r-1}$.

\bigskip
\noindent{\bf 2. General background and
notations.}
\medskip

We fix $X$ and $Y$ as in Theorem 1. If $U$ is an
affine subset of $X$, the set ${\cal O}_X(U)$ of
regular function on $U$ is the coordinate ring
of $U$ (considered as an affine subset in some
affine space). Observe that then the Hilbert
Nullstellensatz implies that, for any set of
regular functions $f,g\in{\cal O}_X(U)$, it
holds:
$$\{p\in U\ |\ f(p)=g(p)=0\}=\emptyset\Rightarrow
\hbox{there exist }u,v\in{\cal O}_X(U)
\hbox{ such that } uf+vg=1\eqno{(4)}$$
(the same is true for an arbitrary number of
functions, but we will not use it).

The fact that $Y$ is a local complete
intersection subscheme of
$X$ of codimension two implies (the reader who is
not familiar with the theory of schemes can take
this as a definition) that any point of $Y$ has
an affine open neighborhood $U\subset X$ such
that ${\cal J}(U)$, the ideal of $Y\cap U$
inside $U$, is generated by two regular functions
$f,g\in{\cal O}_X(U)$ ``without common
components'' i.e. for any regular functions
$u,v\in{\cal O}_X(U)$ it holds:
$$uf=vg\Rightarrow\hbox{there exists }
w\in{\cal O}_X(U)\hbox{ such that } u=wg,\ v=wf.
\eqno{(5)}$$

Instead of regarding vector bundles as locally
free sheaves (as we did in the previous
section), we will consider their geometric
interpretation. Hence, for a vector bundle $E$ of
rank $r$ over an algebraic variety $X$ we will
take an (affine) open covering
$X=\bigcup_{i\in I}U_i$ such that $E_{|U_i}\cong
U_i\times k^r$ (i.e. $E$ {\it trivializes}
on $U_i$). For any $i,j\in I$, elements in
$E_{|U_i\cap U_j}$ can be regarded as elements
in both $U_j\times k^r$ and $U_i\times k^r$, and
the pass from one to another is given by the
multiplication by an $r\times r$ {\it transition
matrix} $Z_{ij}$ of regular functions on $U_i\cap
U_j$ (when $r=1$, we just speak of the {\it transition
functions} of the line bundle). Hence a vector
bundle can be characterized by a collection of
matrices $\{Z_{ij}\}_{i,j\in I}$ subject to the
compatibility condition $Z_{ik}=Z_{ij}Z_{jk}$
(and $Z_{ii}=I_r$, the identity matrix). 

If $U$ is an affine set of $X$, then $Y\cap U$
is also affine, and hence
${\cal O}_Y(Y\cap U)={\cal O}_X(U)/(f,g)$. We
will always use a bar to indicate the classes of
elements (vector bundles, functions,
matrices,...) modulo $Y$ (or any $Y\cap U$).

\bigskip
\noindent{\bf 3. Affine coverings of $X$}
\medskip

We start taking a covering of $Y$ by affine
sets $Y\cap U_i$ (with $i$ varying in a set $I$)
such that: 
\item{(i)} $U_i$ is an affine set of $X$.
\item{(ii)} The vector bundle $L$ trivializes on
$U_i$ and has transition functions $h_{ij}$.
\item{(iii)} ${\cal J}(U_i)$ is generated by
the vanishing of two regular functions $f_i,g_i$
on $U_i$.

In the intersection of two of those open sets,
$U_i,U_j$ we have now two different sets of
generators for the ideal ${\cal J}(U_i\cap
U_j)$, and hence it is possible to find a matrix
$A_{ij}$ (not necessarily unique) satisfying:  
$$\pmatrix{f_i\cr g_i}=
A_{ij}\pmatrix{f_j\cr g_j}=
\pmatrix{a_{ij}&b_{ij}\cr c_{ij}&d_{ij}}
\pmatrix{f_j\cr g_j} \eqno{(6)}$$
where $a_{ij},b_{ij},c_{ij},d_{ij}$ are regular
functions on $U_i\cap U_j$ and $\det A_{ij}$
does not have zeros on $U_i\cap U_j$. Observe
that it could happen that $Y\cap U_i\cap
U_j=\emptyset$. In this case, by (4), we can
find $u_i,v_i,u_j,v_j$ such that
$u_if_i+v_ig_i=1=u_jf_j+v_jg_j$, and thus we can
take $A_{ij}=\pmatrix{f_i&-v_i\cr g_i&u_i}
\pmatrix{u_j&v_j\cr -g_j&f_j}$.

Observe also that the vector bundle $N$
trivializes on $Y\cap U_i$ and has as transition
matrices the restriction $\bar A_{ij}$ of
$A_{ij}$ to $Y\cap U_i\cap U_j$.

Let $s_1,\ldots,s_{r-1}$ be the global sections
generating $\bigwedge^2N\otimes L^*$.
For $t=1,\ldots,r-1$, the section $s_t$ can be
represented locally at each $Y\cap U_i$ by a
regular function $\bar s_{it}$ such that there
are relations
$$\bar s_{it}={\det\bar A_{ij}\over\bar
h_{ij}}\bar s_{jt}\eqno{(7)}$$

Since $\bar s_{i1},\ldots\bar
s_{i,r-1}$ do not vanish simultaneously on
$Y\cap U_i$, we can refine the covering and
assume that there is $t_i\in\{1,\ldots,r-1\}$
such that $\bar s_{it_i}$ does not have zeros
in $Y\cap U_i$. Replacing $U_i$ with its
intersection with $\{s_{it_i}\ne0\}$, we can
assume that $s_{it_i}$ does not have zeros in
$U_i$, i.e. it is a unit in ${\cal O}_X(U_i)$.

\proclaim Lemma 2. With the above notations, it
is possible to choose regular functions
$f_i,g_i$ such that $s_{it_i}=(-1)^{t_i}$. In
particular,
$\det\bar A_{ij}=(-1)^{t_i}{\bar h_{ij}\over\bar
s_{jt_i}}$.

\noindent{\it Proof:} 
We choose as new set of generators of each
${\cal J}(U_i)$ the functions 
$f'_i={f_i\over s_{it_i}}$ and
$g'_i=(-1)^{t_i}g_i$. We obtain a new relation
like (6) with a new matrix
$A'_{ij}$:
$$\pmatrix{f'_i\cr g'_i}=
A'_{ij}\pmatrix{f'_j\cr g'_j}=
\pmatrix{
{s_{jt_j}\over s_{it_i}}a_{ij}&{(-1)^{-t_j}\over
s_{it_i}}b_{ij}\cr
(-1)^{t_i}c_{ij}s_{jt_j}&(-1)^{t_i-t_j}d_{ij} }
\pmatrix{f'_j\cr g'_j}$$ 
from which we get, by (7),
$$(-1)^{t_i}{\bar s_{it}\over\bar s_{it_i}}=
(-1)^{t_j}{\det\bar A'_{ij}\over\bar h_{ij}}
{\bar s_{jt}\over\bar s_{jt_j}}.$$
This shows that, with this new choice of
$f'_i,g'_i$, the sections $s_1,\ldots,s_{r-1}$
can be represented in $Y\cap U_i$ by the
classes of $(-1)^{t_i}{s_{i1}\over
s_{it_i}},\ldots,(-1)^{t_i}{s_{i,r-1}\over
s_{it_i}}$. This implies that we can assume
$s_{it_i}=(-1)^{t_i}$. With this choice, the
last statement is just (7) applied to $t=t_i$.
\qed

\bigskip

We extend now the affine covering to a covering
of the whole $X$. For this, we have to cover
$X\setminus Y$ by new affine open sets $U_i$.
For such a new open set we take $f_i=1,g_i=0$.
Observe that, even if $Y\cap U_i=\emptyset$,
property (5) still holds in a trivial way. 

We also have matrices $A_{ij}$ as in (6)
for any choice of open sets $U_i,U_j$.
Specifically:

--If $Y\cap U_i\ne\emptyset\ne Y\cap U_j$, we do
as in (6).

--If $Y\cap U_i=\emptyset=Y\cap U_j$, we take
$A_{ij}$ to be the identity matrix.

--If $Y\cap U_i\ne\emptyset=Y\cap U_j$, we
take $A_{ij}=\pmatrix{u_j&\hskip -3mm v_j\cr 
-g_j&\hskip -3mm f_j}$, with $u_j,v_j$ such that
$u_jf_j+v_jf_j=1$.

--If $Y\cap U_i=\emptyset\ne Y\cap U_j$, we
take $A_{ij}=\pmatrix{f_i&-v_i\cr g_i&u_i}$,
with $u_i,v_i$ such that $u_if_i+v_if_i=1$.

\proclaim Lemma 3. With the above choices and
notations, it is possible to choose the matrices
$A_{ij}$ such that $\det
A_{ij}=(-1)^{t_i}{h_{ij}\over s_{jt_i}}$.

\noindent{\it Proof.} By Lemma 2, on each
$U_i\cap U_j$ the regular functions $\det
A_{ij}$ and
$(-1)^{t_i}{h_{ij}\over s_{jt_i}}$ coincide
modulo the ideal $(f_i,g_i)$ (this is trivial if
$Y\cap U_i=\emptyset$). We can thus write
$$(-1)^{t_i}{h_{ij}\over s_{jt_i}}=\det
A_{ij}+\varphi_{ij}f_i+\psi_{ij}g_i=
\det A_{ij}
+(\varphi_{ij}a_{ij}+\psi_{ij}c_{ij})f_j
+(\varphi_{ij}b_{ij}+\psi_{ij}d_{ij})g_j$$
for some regular functions $\varphi_{ij},
\psi_{ij}$ on $U_i\cap U_j$. Therefore we can
replace (6) with
$$\pmatrix{f_i\cr g_i}=\pmatrix{
a_{ij}+\psi_{ij}g_j&b_{ij}-\psi_{ij}f_j\cr
c_{ij}-\varphi_{ij}g_j&d_{ij}+\varphi_{ij}f_j}
\pmatrix{f_j\cr g_j} $$
and the new transition matrix $A'_{ij}=\pmatrix{
a_{ij}+\psi_{ij}g_j&b_{ij}-\psi_{ij}f_j\cr
c_{ij}-\varphi_{ij}g_j&d_{ij}+\varphi_{ij}f_j}$
satisfies the wanted property $\det
A'_{ij}=(-1)^{t_i}{h_{ij}\over s_{jt_i}}$.
\qed

\bigskip
\noindent{\bf 4. Constructing the sections.}
\medskip 

We start by fixing a notation that we will use
in the rest of the paper.

\noindent{\bf Notation.} Given the identity
matrix (whose order will be clear any time from
the context), we will denote by $\Delta_t$ the
submatrix obtained by removing its $t$-th row.
Hence, for any matrix $M$, the matrix
$\Delta_tM$ will be the submatrix of
$M$ obtained by removing its $t$-th row.
Similarly, if $\Delta'_t$ is the transpose of
$\Delta_t$, then $M\Delta'_t$ will be the
submatrix of $M$ obtained by removing its $t$-th
column. 

Before constructing the vector bundle $E$ and its
$r-1$ sections $\alpha_1,\ldots,\alpha_{r-1}$,
let us assume that they exist and see the form
they can take. Assume, without loss of
generality, that $E$ trivializes on each $U_i$.
Since
$s_1{\alpha_1}_{|Y}+\ldots+s_{r-1}{\alpha_{r-1}}_{|Y}=0$
and $s_{t_i}$ is represented by $(-1)^{t_i}$ on
$Y\cap U_i$, this means that, on the points of $Y$,
$\alpha_{t_i}$ depends on
$\alpha_1,\ldots,\hat\alpha_{t_i}\ldots,\alpha_{r-1}$
(we use the standard notation of a hat to
indicate that a term is removed). Since the rank
of $\alpha_1,\ldots,\alpha_{r-1}$ is
$r-2$ on $Y$ and $r-1$ outside $Y$, it follows that
$\alpha_1,\ldots,\hat\alpha_{t_i}\ldots,\alpha_{r-1}$
are linearly independent on $U_i$. Extending
them to a basis of $E_{|U_i}$, it is then
possible to represent
$\alpha_1,\ldots,\alpha_{r-1}$ on $U_i$, in
terms of this basis, as the columns of an
$r\times(r-2)$ matrix
$M_i=\Delta_{t_i}T_i$, where
$$T_i=\pmatrix{
1&0&\ldots&\alpha_{i1}&\ldots&0\cr
0&1&\ldots&\alpha_{i2}&\ldots&0\cr
\vdots&&\ddots&\vdots&&\vdots\cr
0&0&\ldots&\alpha_{it_i}&\ldots&0\cr
\vdots&&&\vdots&\ddots&\vdots\cr
0&0&\ldots&\alpha_{i,r-1}&\ldots&1\cr
0&0&\ldots&\alpha_{ir}&\ldots&0\cr
0&0&\ldots&\alpha_{i,r+1}&\ldots&0}$$
Since $Y\cap U_i$ must the determinantal variety
defined by the maximal minors of $M_i$, it follows
that $\alpha_{ir},\alpha_{i,r+1}$ generate
${\cal J}(U_i)$. Hence, changing the last two
rows of $T_i$ by a suitable linear combination of
them, we can assume
$\alpha_{ir}=f_i,\ \alpha_{i,r+1}=g_i$. 

On the other hand, the equation $s_1{\alpha_1}_{|Y}+\ldots+s_{r-1}{\alpha_{r-1}}_{|Y}=0$
implies that the entries of
$M_i\pmatrix{s_{i1}\cr\vdots\cr
s_{i,r-1}}$ are a linear combination of
$f_i,g_i$. Hence, after adding to each of the
first $r-2$ rows of $M_i$ a linear combination of
the last two, we can take 
$$T_i=\pmatrix{T'_i\cr T''_i}$$
with 
$$T'_i=\pmatrix{
1&0&\ldots&-(-1)^{t_i}s_{i1}&\ldots&0\cr
0&1&\ldots&-(-1)^{t_i}s_{i2}&\ldots&0\cr
\vdots&&\ddots&\vdots&&\vdots\cr
0&0&\ldots&1&\ldots&0\cr
\vdots&&&\vdots&\ddots&\vdots\cr
0&0&\ldots&-(-1)^{t_i}s_{i,r-1}&\ldots&1}\eqno{(8)}$$
and
$$T''_i=\pmatrix{
0&0&\ldots&f_i&\ldots&0\cr
0&0&\ldots&g_i&\ldots&0}.\eqno{(9)}$$

We will thus define
$$M_i=\pmatrix{\Delta_{t_i}T'_i\cr T''_i}
\eqno{(10)}$$
with $T'_i$ and $T''_i$ as in (8) and (9). We
have the following easy equalities, which we
will use frequently:

$$\Delta_{t_i}T'_i\Delta'_{t_i}=I_{r-2}
\eqno{(11)}$$

$$T''_i\Delta'_{t_i}=
\pmatrix{0&\ldots&0\cr0&\ldots&0}
\eqno{(12)}$$

$$\Delta_{t_i}T'_i\pmatrix{s_{i1}\cr\vdots\cr
s_{i,r-1}}=
\pmatrix{0\cr\vdots\cr0}
\eqno{(13)}$$ 

$$T''_i\pmatrix{s_{i1}\cr\vdots\cr s_{i,r-1}}=
(-1)^{t_i}\pmatrix{f_i\cr g_i}.
\eqno{(14)}$$ 

Since we want the columns of $M_i$ to represent
the sections $\alpha_1,\ldots,\alpha_{r-1}$ of a
vector bundle $E$, we need to find the
transition matrices relating $M_i$ to $M_j$.
The next result provides a first condition to
find them.

\proclaim Lemma 4. For a covering and choices as
in Lemma 3, if for each $i\in I$ we take
$M_i$ as in (10), then an $r\times r$ matrix
$Z_{ij}=\pmatrix{P_{ij}&Q_{ij}\cr
R_{ij}&S_{ij}\cr}$ satisfies the equality
$M_i=Z_{ij}M_j$ if and only if the following
equalities hold:
\item{(i)} $P_{ij}=\Delta_{t_i}T'_i\Delta'_{t_j}$
\item{(ii)} $R_{ij}=T''_i\Delta'_{t_j}$
\item{(iii)} $Q_{ij}\pmatrix{f_j\cr
g_j}=(-1)^{t_j}
\Delta_{t_i}T'_i\pmatrix{s_{j1}\cr\vdots\cr
s_{j,r-1}}$
\item{(iv)} $S_{ij}\pmatrix{f_j\cr g_j}
=(-1)^{t_j}s_{jt_i}\pmatrix{f_i\cr g_i}$, i.e.
$S_{ij}=(-1)^{t_j}s_{jt_i}A_{ij}$, with
$A_{ij}$ as in (6). 
\smallskip\noindent Moreover, such a matrix
always exists and, when taking $A_{ij}$ as in
Lemma 3, it follows $\det
S_{ij}=(-1)^{t_i}s_{jt_i}h_{ij}$ and $\det
Z_{ij}=h_{ij}$.

\noindent{\it Proof.} We have to find
the solutions of 
$$\left\{\eqalign{\Delta_{t_i}T'_i=P_{ij}\Delta_{t_j}T'_j+Q_{ij}T''_j\cr
T''_i=R_{ij}\Delta_{t_j}T'_j+S_{ij}T''_j}\right.\eqno{(15)}$$
Multiplying by $\Delta'_{t_j}$ to the right the
two equations in (15) (i.e. removing the
$t_j$-th columns of all the terms), we get from
(11) and (12) the equalities (i) and (ii). It
remains to characterize when (15) holds for the
$t_j$-th column of each term. To see this,
since $s_{jt_j}=(-1)^{t_j}$, it is equivalent to
consider the product of the two equalities of
(15) with
$\pmatrix{s_{j1}\cr\vdots\cr s_{j,r-1}}$,
which together with (13) and (14) yield exactly
the equalities (iii) and (iv). 

The entries of
$\Delta_{t_i}T'_i\pmatrix{s_{j1}\cr\vdots\cr
s_{j,r-1}}$ are 
$s_{jt}-(-1)^{t_i}s_{it}s_{jt_i}$, with
$t=1,\ldots,\hat t_i,\ldots,r-1$. Recalling
from Lemma 3 that
$\det A_{ij}=(-1)^{t_i}{h_{ij}\over s_{jt_i}}$,
equality (7) reads $\bar s_{jt}-(-1)^{t_i}\bar
s_{it}\bar s_{jt_i}=\bar0$. Hence the entries of
$\Delta_{t_i}T'_i\pmatrix{s_{j1}\cr\vdots\cr
s_{j,r-1}}$ are in the ideal $(f_j,g_j)$
defining $Y\cap U_i\cap U_j$, and the same holds
clearly for the entries of
$s_{jt_i}\pmatrix{f_i\cr g_i}$. Therefore,
equalities (iii) and (iv) have solutions
$Q_{ij},S_{ij}$, and thus there exists some
$Z_{ij}$ such that
$M_i=Z_{ij}M_j$.

For the last equality in the statement, we
deduce from the equations (15) multiplied to the
right by $\Delta'_{t_i}$, and using (11) and
(12), the equality
$$\pmatrix{P_{ij}&Q_{ij}\cr R_{ij}&S_{ij}\cr}
\pmatrix{\Delta_{t_j}T'_j\Delta'_{t_i}&0\cr
T''_j\Delta'_{t_i}&I_2}=
\pmatrix{I_{r-2}&Q_{ij}\cr 0&S_{ij}}.
$$
Hence, observing that
$\det(\Delta_{t_j}T'_j\Delta'_{t_i})=(-1)^{t_i}s_{jt_i}$,
we obtain
$(-1)^{t_i}s_{jt_i}\det Z_{ij}=\det S_{ij}$.
Since $S_{ij}=(-1)^{t_i}s_{jt_i}A_{ij}$ and
$\det A_{ij}=(-1)^{t_j}{h_{ij}\over s_{jt_i}}$
after Lemma 3, we thus have $\det
S_{ij}=(-1)^{t_i}s_{jt_i}h_{ij}$ and therefore
$\det Z_{ij}=h_{ij}$.
\qed

\noindent{\bf Remark 5.} It is only now that
one can understand the reason of introducing
the sign $(-1)^{t_i}$ in Lemma 2. Observe first
that it was not a misprint to write
$\det(\Delta_{t_j}T'_j\Delta'_{t_i})=(-1)^{t_i}s_{jt_i}$
at the end of the proof of Lemma 4, in the
sense that it is indeed $(-1)^{t_i}$ instead of
$(-1)^{t_j}$ (which is the sign appearing in the
entries of the matrix $T'_j$). If we had not
included that sign in Lemma 2, we would have
obtained now $\det Z_{ij}=(-1)^{t_i+t_j}h_{ij}$
in Lemma 4. This would not have been a disaster,
since the functions $(-1)^{t_i+t_j}h_{ij}$ are
also transition functions of $L$. Anyway, we
thought it was more elegant and clearer not to
work simultaneously with two different sets of
transition functions of the same line bundle.

\proclaim Lemma 6. For a matrix $Z_{ij}$ as in
Lemma 4, the following equalities hold:
\item{(i)}
$(g_i,-f_i)S_{ij}=(-1)^{t_i+t_j}h_{ij}(g_j,-f_j)$.
\item{(ii)} $R_{ij}=\pmatrix{f_i\cr
g_i}(\delta_{ij1}\ldots\hat\delta_{ijt_j}\ldots\delta_{ij,r-1})$,
with $\delta_{ijt}=0$, for all $t\ne t_i$ and
$\delta_{ijt_i}=1$; in particular,
$(g_i,-f_i)R_{ij}=(0\ldots 0)$.
\item{(iii)}
$(0\ldots0,g_i,-f_i)Z_{ij}=(-1)^{t_i+t_j}h_{ij}(0\ldots
0,g_j,-f_j)$.
\item{(iv)} $(0\ldots0,g_i,-f_i)(Z_{ij}Z_{jk}-Z_{ik})=
(0\ldots0)$.

\noindent{\it Proof.} The equality
$S_{ij}\pmatrix{f_j\cr
g_j}=(-1)^{t_j}s_{jt_i}\pmatrix{f_i\cr g_i}$ of
Lemma 4 is equivalent, multiplying to
the left by Adj$S_{ij}$ and using $\det
S_{ij}=(-1)^{t_i}s_{jt_i}h_{ij}$, to
$\pmatrix{f_j\cr g_j}={(-1)^{t_i+t_j}\over
h_{ij}}{\rm Adj}S_{ij}\pmatrix{f_i\cr g_i}$,
which is in turn equivalent to (i). Part (ii) is
obvious, since $R_{ij}=T''_i\Delta'_{t_j}$.
Part (iii) follows from (i) and (ii). Finally,
part (iv) is a consequence of (iii), having in
mind, by Lemma 4, that the $h_{ij}$ are the
transition functions of the line bundle $L$ and
therefore
$h_{ij}h_{jk}=h_{ik}$. 
\qed

\proclaim Corollary 7. If the matrices
$\{Z_{ij}\}_{i,j,\in I}$ are chosen as in Lemma
4, then for any $i,j,k\in I$ there exist
regular functions
$\beta_{ijk1},\ldots,\beta_{ijk,r-1}$ on
$U_i\cap U_j\cap U_k$ such that
$Z_{ik}-Z_{ij}Z_{jk}=(0\ B_{ijk})$, with
$$B_{ijk}=\pmatrix{Q_{ik}-P_{ij}Q_{jk}-Q_{ij}S_{jk}\cr
\cr S_{ik}-R_{ij}Q_{jk}-S_{ij}S_{jk}}
=\pmatrix{\beta_{ijk1}\cr\vdots\cr
\hat\beta_{ijkt_i}\cr\vdots\cr
\beta_{ijk,r-1}\cr
\beta_{ijkt_i}f_i\cr\beta_{ijkt_i}g_i}(g_k,-f_k)$$ 

\noindent{\it Proof.} Write $Z_{ik}-Z_{ij}Z_{jk}=(B'_{ijk}\
B''_{ijk})$. The equality
$(Z_{ik}-Z_{ij}Z_{jk})M_k=0$ is equivalent to
$B'_{ijk}\Delta_{t_k}T'_k+B''_{ijk}T''_k=0$,
so it follows, multiplying this equality to the
right by $\Delta'_{t_k}$ and applying (11) and
(12), that $B'_{ijk}=0$. Hence 
$B''_{ijk}T''_k=0$, i.e., by the definition (9) of $T''_k$,
$B''_{ijk}\pmatrix{f_k\cr fg_k}=0$. It follows from (5)
that there exist regular functions
$\beta_{ijk1},\ldots,\hat\beta_{ijkt_i},\ldots,
\beta_{ijk,r+1}$ such that
$B_{ijk}=\pmatrix{
\beta_{ijk1}\cr\vdots\cr\hat\beta_{ijkt_i}\cr\vdots\cr
\beta_{ijk,r+1}}(g_k,-f_k)$. On the other
hand, applying now Lemma 6(iv), we get $(g_i,\
-f_i)\pmatrix{\beta_{ijkr}\cr\beta_{ijk,r+1}}=0$,
from which the lemma follows by applying
(5) again.
\qed

\noindent{\bf Remark 8.} If we want the matrices
$Z_{ij}$ to be the transition matrices of a
vector bundle $E$, we need to find a good choice
of $Q_{ij},S_{ij}$ such that
$\beta_{ijk1},\ldots,\beta_{ijk,r-1}$ are all
zero. Observe that another choice of
$Q'_{ij}$ and $S'_{ij}$ satisfies 
conditions (iii) and (iv) of Lemma 4 if and
only if we have respectively
$(Q'_{ij}-Q_{ij})\pmatrix{f_j\cr
g_j}=\pmatrix{0\cr\vdots\cr0}$ and
$(S'_{ij}-S_{ij})\pmatrix{f_j\cr
g_j}=\pmatrix{0\cr0}$. Moreover, using Lemma
6(i), we would also have $(g_i,\
-f_i)(S'_{ij}-S_{ij})=(0\ 0)$. Hence the same
reasoning as in the proof of Corollary 7 shows
that the above conditions are equivalent to the
existence of regular functions
$x_{ij1},\ldots,x_{ij,r-1}$ such that
$Q'_{ij}=Q_{ij}+\pmatrix{x_{ij1}\cr\vdots\cr\hat
x_{ijt_i}\cr\vdots\cr x_{ij,r-1}}(g_j,\ -f_j)$
and
$S'_{ij}=S_{ij}+x_{ijt_i}\pmatrix{f_i\cr
g_i}(g_j,\ -f_j)$. The goal of the next section
will be to see that there is essentially one way
of choosing the functions
$x_{ij1},\ldots,x_{ij,r-1}$ on each $U_i\cap U_j$. We will
then see how the a priori strange choice of subindices
makes perfectly sense.

\bigskip
\noindent{\bf 5. Constructing the vector bundle.}
\medskip

We finally find under which conditions the
matrices $Z_{ij}$ are transition matrices of a
vector bundle. We start with a technical lemma
that will be very useful in the sequel:

\proclaim Lemma 9. With the definitions of
the previous section, for any vector
$u=\pmatrix{u_1\cr\vdots\cr u_{r-1}}$, we have
${T'_j}^{-1}u=\Delta'_{t_j}\pmatrix{
u_1\cr\vdots\cr\hat u_{t_j}\cr\vdots\cr
u_{r-1}}+(-1)^{t_j}u_{t_j}\pmatrix{s_{j1}\cr\vdots\cr
s_{j,r-1}}$. Hence, if we
define $u'=T'_i{T'_j}^{-1}u$, then:
\item{(i)} $\Delta_{t_i}u'=
P_{ij}\pmatrix{u_1\cr\vdots\cr\hat u_{t_j}\cr\vdots\cr u_{r-1}}+
(-1)^{t_j}u_{t_j}\Delta_{t_i}T'_i\pmatrix{s_{j1}\cr\vdots\cr
s_{j,r-1}}$.
\item{(ii)} The $t_i$-th row of $u'$ is $
(\delta_{ij1}\ldots\hat\delta_{ijt_j}\ldots\delta_{ij,r-1})
\pmatrix{u_1\cr\vdots\cr\hat u_{t_j}\cr\vdots\cr u_{r-1}}
+(-1)^{t_j}s_{jt_i}u_{t_j}$.

\noindent{\it Proof.} For the first equality,
observe first that we can write
$u=\Delta'_{t_j}\pmatrix{ u_1\cr\vdots\cr\hat
u_{t_j}\cr\vdots\cr
u_{r-1}}+u_{t_j}\pmatrix{0\cr\vdots\cr1\cr\vdots\cr
0}$ (the first summand is nothing but $u$ with
the $t_j$-th row replaced with $0$). Then the
wanted equality follows because
${T'_j}^{-1}=\pmatrix{
1&0&\ldots&(-1)^{t_j}s_{j1}&\ldots&0\cr
0&1&\ldots&(-1)^{t_j}s_{j2}&\ldots&0\cr
\vdots&&\ddots&\vdots&&\vdots\cr
0&0&\ldots&1&\ldots&0\cr
\vdots&&&\vdots&\ddots&\vdots\cr
0&0&\ldots&(-1)^{t_j}s_{j,r-1}&\ldots&1}$
and then
${T'_j}^{-1}\Delta'_{t_j}=\Delta'_{t_j}$
(observe also that $(-1)^{t_j}s_{jt_j}=1$) .

Now (i) and (ii) are easy consequences of the
first equality: for (i) it is enough to recall
from Lemma 4 that
$P_{ij}=\Delta_{t_i}T'_i\Delta'_{t_j}$, while
for (ii) it suffices to observe that
the $t_i$-th row of $\Delta'_{t_j}$ is $
(\delta_{ij1}\ldots\hat\delta_{ijt_j}\ldots\delta_{ij,r-1})
$.
\qed

\proclaim Proposition 10. For a choice of
matrices $Z_{ij}$ as in Lemma 4, let
$Z'_{ij}=\pmatrix{P_{ij}&Q'_{ij}\cr
R_{ij}&S'_{ij}}$ with $Q'_{ij}=Q_{ij}+\pmatrix{x_{ij1}\cr\vdots\cr\hat
x_{ijt_i}\cr\vdots\cr x_{ij,r-1}}(g_j,\ -f_j)$
and
$S'_{ij}=S_{ij}+x_{ijt_i}\pmatrix{f_i\cr
g_i}(g_j,\ -f_j)$. Then
$Z'_{ik}-Z'_{ij}Z'_{jk}=0$ if and only if 
$$\eqalign{
(-1)^{t_k}{T'_i}^{-1}\pmatrix{\beta_{ijk1}\cr\vdots\cr\beta_{ijk,r-1}}=
(-1)^{t_k}{T'_j}^{-1}\pmatrix{x_{jk1}\cr\vdots\cr
x_{jk,r-1}}-\cr
&\hskip -5cm
-(-1)^{t_k}{T'_i}^{-1}\pmatrix{x_{ik1}\cr\vdots\cr
x_{ik,r-1}}+
(-1)^{t_j}h_{jk}{T'_i}^{-1}\pmatrix{x_{ij1}\cr\vdots\cr
x_{ij,r-1}}.}
\eqno{(16)}$$

\noindent{\it Proof.} Multiplying to the
left by $(-1)^{t_k}T'_i$, equation (16) in the
statement is equivalent to
$$\pmatrix{\beta_{ijk1}\cr\vdots\cr\beta_{ijk,r-1}}
-{T'_i}{T'_j}^{-1}\pmatrix{x_{jk1}\cr\vdots\cr
x_{jk,r-1}}
+\pmatrix{x_{ik1}\cr\vdots\cr x_{ik,r-1}}
-(-1)^{t_j+t_k}h_{jk}\pmatrix{x_{ij1}\cr\vdots\cr
x_{ij,r-1}}
=\pmatrix{0\cr\vdots\cr0}.$$ 
Looking separately to the $t_i$-th row and the
others, the above equality is equivalent, by
Lemma 9, to the vanishing of
$$\Lambda_{ijk}:=
\pmatrix{\beta_{ijk1}\cr\vdots\cr\hat\beta_{ijkt_i}
\cr\vdots\cr\beta_{ijk,r-1}}
-P_{ij}\pmatrix{x_{jk1}\cr\vdots\cr\hat
x_{jkt_j}\cr\vdots\cr x_{jk,r-1}}-\hskip 6cm$$
$$\hskip 2cm
-(-1)^{t_j}x_{jkt_j}\Delta_{t_i}T'_i\pmatrix{s_{j1}\cr\vdots\cr
s_{j,r-1}}
+\pmatrix{x_{ik1}\cr\vdots\cr
\hat x_{ikt_i}\cr\vdots\cr
x_{ik,r-1}}
-(-1)^{t_j+t_k}h_{jk}\pmatrix{x_{ij1}\cr\vdots\cr
\hat x_{ijt_i}\cr\vdots\cr x_{ij,r-1}}
$$
and
$$\lambda_{ijk}:=\beta_{ijkt_i}
-(\delta_{ij1}\ldots\hat\delta_{ijt_j}
\ldots\delta_{ij,r-1})\pmatrix{x_{jk1}
\cr\vdots\cr\hat x_{jkt_j}\cr\vdots\cr
x_{jk,r-1}}-
\hskip 4cm$$
$$\hskip 5cm 
-(-1)^{t_j}s_{jt_i}x_{jkt_j}+x_{ikt_i}
-(-1)^{t_j+t_k}h_{jk}x_{ijt_i}.
$$

On the other hand, the condition
$Z'_{ik}-Z'_{ij}Z'_{jk}=0$ is equivalent, by
Corollary 7, to the vanishing of 
$Q'_{ik}-P_{ij}Q'_{jk}-Q'_{ij}S'_{jk}$ and
$S'_{ik}-R_{ij}Q'_{jk}-S'_{ij}S'_{jk}$. A 
straightforward calculation (using Lemmas 4 and
6 and Corollary 7) shows that

$$Q'_{ik}-P_{ij}Q'_{jk}-Q'_{ij}S'_{jk}=\Lambda_{ijk}(g_k,-f_k)$$
and
$$S'_{ik}-R_{ij}Q'_{jk}-S'_{ij}S'_{jk}=
\lambda_{ijk}\pmatrix{f_i\cr g_i}
(g_k,\ -f_k)$$
so that the lemma follows at once.
\qed

\noindent{\bf Remark 11.} Equality (16) means
that the $(r-1)$-uples
$(-1)^{t_k}{T'_i}^{-1}\pmatrix{\beta_{ijk1}\cr\vdots\cr\beta_{ijk,r-1}}$
represents a $2$-coboundary in the {\v C}ech
cohomology of $(r-1)L^*$ with respect to the
covering $\{U_i\}_{i\in I}$ (multiplication by
$h_{jk}$ in the last summand is needed in order
to have all the $(r-1)$-uples defined in the
trivialization of $(r-1)L^*$ in $U_k$). Recall
(see [8] III-Theorem 4.5) that the cohomology
of coherent sheaves (and in particular of
vector bundles) is isomorphic to the {\v C}ech
cohomology of any affine cover. Hence the
matrices $Z'_{ij}$ will be the transition
matrices of a vector bundle as soon as we see
that the $(r-1)$-uples
$(-1)^{t_k}{T'_i}^{-1}\pmatrix{\beta_{ijk1}\cr\vdots\cr\beta_{ijk,r-1}}$
represent a $2$-cocycle, since we are assuming
$H^2(X,(r-1)L^*)=0$. This is what we are going
to do next.

\proclaim Proposition 12. The set of
$(r-1)$-uples
$(-1)^{t_k}{T'_i}^{-1}\pmatrix{\beta_{ijk1}\cr\vdots\cr\beta_{ijk,r-1}}$
defines a $2$-cocycle of the vector
bundle $(r-1)L^*$.

\noindent{\it Proof.} We need to show that for
each $i,j,k,l\in I$, it follows:
$$(-1)^{t_l}{T'_j}^{-1}\pmatrix{\beta_{jkl1}\cr\vdots\cr\beta_{jkl,r-1}}-
(-1)^{t_l}{T'_i}^{-1}\pmatrix{\beta_{ikl1}\cr\vdots\cr\beta_{ikl,r-1}}+
\hskip 4cm$$
$$\hskip 3cm
+(-1)^{t_l}{T'_i}^{-1}\pmatrix{\beta_{ijl1}\cr\vdots\cr\beta_{ijl,r-1}}-
(-1)^{t_k}h_{kl}{T'_i}^{-1}\pmatrix{\beta_{ijk1}\cr\vdots\cr\beta_{ijk,r-1}}=
\pmatrix{0\cr\vdots\cr0}.$$
As in the proof of Proposition 10, multiplying
to the left by $(-1)^{t_l}T'_i$ and applying
Lemma 9, the above equality is equivalent to the
vanishing of
$$\Delta_{ijkl}:=P_{ij}\pmatrix{\beta_{jkl1}\cr\vdots\cr
\hat\beta_{jklt_j}\cr\vdots\cr
\beta_{jkl,r-1}\cr}
+(-1)^{t_j}\beta_{jklt_j}\Delta_{t_i}T'_i\pmatrix{s_{j1}\cr\vdots\cr
s_{j,r-1}}-\hskip 4cm$$
$$\hskip 3cm-\pmatrix{\beta_{ikl1}\cr\vdots\cr
\hat\beta_{iklt_i}\cr\vdots\cr
\beta_{ikl,r-1}\cr}+\pmatrix{\beta_{ijl1}\cr\vdots
\cr\hat\beta_{ijkt_i}\cr\vdots\cr
\beta_{ijl,r-1}\cr}
-(-1)^{t_k+t_l}h_{kl}
\pmatrix{\beta_{ijk1}\cr
\vdots\cr\hat\beta_{ijkt_i}\cr\vdots\cr
\beta_{ijk,r-1}\cr}
$$
and 
$$\lambda_{ijkl}:=(\delta_{ij1}\ldots\hat\delta_{ijt_j}\ldots\delta_{ij,r-1})
\pmatrix{\beta_{jkl1}\cr\vdots\cr\hat\beta_{jklt_j}\cr\vdots\cr
\beta_{jkl,r-1}\cr}+\hskip 4cm$$
$$\hskip 3cm+(-1)^{t_j}\beta_{jklt_j}s_{jt_i}
-\beta_{iklt_i}+\beta_{ijlt_i}
-(-1)^{t_k+t_l}h_{kl}\beta_{ijkt_i}.$$

To prove those equalities, we use the equality
$$Z_{ij}(Z_{jl}-Z_{jk}Z_{kl})
-(Z_{il}-Z_{ik}Z_{kl})
+(Z_{il}-Z_{ij}Z_{jl})
-(Z_{ik}-Z_{ij}Z_{jk})Z_{kl}=0.$$ 
Using Corollary 7 to split the above equality in
two blocks --the one of the first $r-2$ rows and
the one of the last $2$ rows-- and applying then
Lemma 4, we get that the equality is equivalent
to the vanishing of the matrices
$\Lambda_{ijkl}(g_l,-f_l)$ and
$\lambda_{ijkl}\pmatrix{f_i\cr
g_i}(g_l,-f_l)$
which proves the proposition.
\qed

\noindent{\bf Remark 13.} Although I did not
check it, it is natural to expect that the
map $\psi$ in (3) assigns to the morphism
defined by $s_1,\ldots,s_{r-1}$ the cocycle of
Proposition 12.

We finally prove the uniqueness statement.

\proclaim Proposition 14. Assume $E$ is a
vector bundle on $X$ satisfying conditions
(i) and (ii) in Theorem 1. If $H^1(X,L^*)=0$,
then any other vector bundle $E'$ satisfying
the same conditions is isomorphic to $E$.

\noindent{\it Proof.} Assume that the
transition matrices of $E$ and $E'$ are (see
Lemma 4) respectively
$$Z_{ij}=\pmatrix{P_{ij}&Q_{ij}\cr
R_{ij}&S_{ij}}$$ and
$$Z'_{ij}=\pmatrix{P_{ij}&Q'_{ij}\cr
R_{ij}&S'_{ij}}$$ 
with (see Remark 8)
$$Q'_{ij}=Q_{ij}+\pmatrix{x_{ij1}\cr\vdots\cr\hat
x_{ijt_i}\cr\vdots\cr x_{ij,r-1}}(g_j,\ -f_j)
\eqno{(17)}$$
and
$$S'_{ij}=S_{ij}+x_{ijt_i}\pmatrix{f_i\cr
g_i}(g_j,\ -f_j).
\eqno{(18)}$$
By Proposition 10, we have
$$(-1)^{t_k}{T'_j}^{-1}\pmatrix{x_{jk1}\cr\vdots\cr
x_{jk,r-1}}-
(-1)^{t_k}{T'_i}^{-1}\pmatrix{x_{ik1}\cr\vdots\cr
x_{ik,r-1}}+
(-1)^{t_j}h_{jk}{T'_i}^{-1}\pmatrix{x_{ij1}\cr\vdots\cr
x_{ij,r-1}}=
\pmatrix{0\cr\vdots\cr0}$$
i.e. the $(r-1)$-uples
$(-1)^{t_j}{T'_i}^{-1}\pmatrix{x_{ij1}\cr\vdots\cr
x_{ij,r-1}}$ define a
$1$-cocycle in $(r-1)L^*$. Since $H^1(X,L^*)=0$, this cocycle is the
coboundary of a $0$-chain defined by
$(r-1)$-uples that we write in the form
$(-1)^{t_i}{T'_i}^{-1}\pmatrix{y_{i1}\cr\vdots\cr
x_{i,r-1}}$. This means
$$(-1)^{t_j}{T'_i}^{-1}\pmatrix{x_{ij1}\cr\vdots\cr
x_{ij,r-1}}
=(-1)^{t_j}{T'_j}^{-1}\pmatrix{y_{j1}\cr\vdots\cr
y_{j,r-1}}-
(-1)^{t_i}h_{ij}{T'_i}^{-1}\pmatrix{y_{i1}\cr\vdots\cr
y_{i,r-1}}.$$ Multiplying as usual the above
relation to the left by $(-1)^{t_j}T'_i$ and
applying Lemma 9 we get that this equality is
equivalent to the vanishing of 
$$\Lambda_{ij}:=
\pmatrix{x_{ij1}\cr\vdots\cr\hat
x_{ijt_i}\cr\vdots\cr x_{ij,r-1}}
-P_{ij}\pmatrix{y_{j1}\cr\vdots\cr\hat
y_{jt_j}\cr\vdots\cr y_{j,r-1}}
-(-1)^{t_j}y_{jt_j}\Delta_{t_i}T'_i\pmatrix{s_{j1}\cr\vdots\cr
s_{j,r-1}}
+(-1)^{t_i+t_j}h_{ij}\pmatrix{y_{i1}\cr\vdots\cr\hat
y_{it_i}\cr\vdots\cr y_{i,r-1}}
$$
and 
$$\lambda_{ij}:=x_{ijt_i}
-(\delta_{ij1}\ldots\hat\delta_{ijt_j}\ldots\delta_{ij,r-1})
\pmatrix{y_{j1}\cr\vdots\cr\hat
y_{yt_j}\cr\vdots\cr
y_{j,r-1}}
-(-1)^{t_j}s_{jt_i}y_{jt_j}
+(-1)^{t_i+t_j}h_{ij}y_{it_i}.
$$
We consider the matrix
$$N_i=\pmatrix{I_{r-2}&N'_i\cr0&I_2+N''_i}$$
where
$$N'_i=\pmatrix{y_{i1}\cr\vdots\cr\hat y_{it_i}\cr\vdots\cr
y_{i,r-1}}(g_i,\ -f_i)\eqno{(19)}$$ 
and 
$$N''_i=y_{it_i}\pmatrix{f_i\cr g_i}(g_i,\
-f_i).\eqno{(20)}$$ 
We define, for each $i\in I$, the automorphism
of the trivial vector bundle $U_i\times k^r$
consisting of the multiplication by $N_i$
(observe that $\det N_i=1$). The result will
be proved if we can patch all these
automorphism in order to get an isomorphism
between $E$ and $E'$. For this, we need to
check the equality $Z_{ij}N_j=N_iZ'_{ij}$.
Splitting this equality in four blocks, it
becomes equivalent to two tautologies (using
Lemma 6(ii)) and the two equalities:
$$P_{ij}N'_j+Q_{ij}+Q_{ij}N''_j=Q'_{ij}+N'_iS'_{ij}$$
and
$$R_{ij}N'_j+S_{ij}+S_{ij}N''_j=S'_{ij}+N''_iS'_{ij}$$
Using (17), (18), (19), (20) and Lemmas 4 and 6,
these two equalities become respectively
equivalent to the vanishing of 
$\Lambda_{ij}(g_j,\ -f_j)$ and 
$\lambda_{ij}\pmatrix{f_i\cr g_i}(g_j,\ -f_j)$,
which completes the proof.
\qed

\noindent{\bf Remark 15.} It is not by chance
that $N_i$ takes the aspect obtained in the
previous proof. It can be easily proved that
this is the aspect that should take any matrix
satisfying $N_iM_i=M_i$ and $\det N_i=1$. In
other words, $N_i$ preserves the local
expression of the sections
$\alpha_1,\ldots,\alpha_{r-1}$ and the
determinant of the transition matrix. This
means that the isomorphism that we found
preserves also the sections
$\alpha_1,\ldots,\alpha_{r-1}$. 

\bigskip
\noindent{\bf 6. Final remarks.}
\medskip

The natural question when trying to
generalize Hartshorne-Serre construction to
higher codimension is:

\noindent{\bf Question 16.} Given a local
complete intersection subscheme $Y$ of
codimension $s$ of a smooth variety $X$, when
is it possible to describe $Y$ as the
dependency locus of $r-s+1$ sections of a rank
$r$ vector bundle $E$ over $X$? When is it
possible to take $r=s$?

If one tries to imitate the technique explained
in section 1, one can regard $Y$ as the
degeneracy locus of a map $V\otimes{\cal
O}_X\to E$, where $V$ is a vector space of
dimension $r-s+1$. The Eagon-Northcott complex
associated to this map produces a long exact
sequence, analog to (1),
$$0\to S^{s-1}V\otimes{\cal O}_X\to
S^{s-2}V\otimes
E\to\ldots\to\bigwedge^{s-1}E\to{\cal J}\otimes
L\to  0\eqno{(21)}$$
where $L=\bigwedge^rE$. Dualizing (21), using
the isomorphism ${\cal E}xt^{s-1}({\cal
J},{\cal O}_X)\cong\bigwedge^sN$ we get an
epimorphism
$$S^{s-1}V^*\otimes{\cal
O}_X\to\bigwedge^sN\otimes L^*.\eqno{(22)}$$
When trying to obtain (22) from (21), as in
section 1 we get that the surjection provides an
element of ${\rm Ext}^{s-1}({\cal J\otimes
L},S^{s-1}V\otimes{\cal O}_X)$. An element
there represents the class of a long exact
sequence of length $s-1$ starting and
finishing as (21), but if $s>2$ the equivalence
classes of these extensions are difficult to
deal with, and it does not look easy to decide
when there is some equivalence class
corresponding to an Eagon-Northcott complex
like (21).

Unfortunately, our construction does not seem
to give a hint to answer Question 16 when $s>2$
neither. Even when $r=s$ (i.e. when we want $Y$
to be the zero locus of a section of a vector
bundle of rank $s$), our construction seems to
suggest that everything could works as soon as
$\bigwedge^2N$ is extendable to $X$, but this is
a very strong condition (for example, if $s=3$
this is essentially equivalent to say that $N$
itself is extendable, which is precisely what we
want to prove). 

Observe also that, in the codimension two case,
Hartshorne-Serre correspondence is saying
(except for the cohomological condition on
$L$) that a local complete intersection subscheme
is the zero locus of a section of as rank two
vector bundle if and only if the Chern classes
of the normal bundle $N$ extend to the ambient
variety (the extendability of the second Chern
class always holds by the self-intersection
formula). However, in higher codimension,
although this condition is clearly necessary
(since $N$ itself has to extend to the ambient
variety) is not at all sufficient (for instance,
most of the elliptic curves in ${\Bbb P}^4$
will provide a counterexample). Hence some
extra condition is needed.

I finally want to mention that we expect that
our construction could be generalized to other
context different from algebraic varieties
over an algebraically closed field. For
instance, property (4) still holds in the
context of real varieties (algebraic or not):
it is enough to take $u={f\over f^2+g^2+1}$
and $v={g\over f^2+g^2+1}$ (I thank Marco
Castrill\'on for suggesting me this idea);
hence the whole construction seems to work in
this new context.

\bigskip
\noindent{\bf References.}
\medskip

\item{[1]} A.B. Altman, S.L. Kleiman, {\it
Introduction to Grothendieck duality theory},
Springer LNM 146, 1970.

\item{[2]} E. Arrondo, {\it La correspondencia
de Serre hecha a mano}, in {\it Homenaje al
profesor Outerelo}, Contribuciones Matem\'aticas,
Editorial Complutense (2004), 61-72.

\item{[3]} W. Barth, A. Van de Ven, {\it A
decomposability criterion for algebraic
2-bundles on projective spaces}, Invent.
Math. 25 (1974), 91-106.

\item{[4]} D. Eisenbud, S. Popescu, C.
Walter, {\it Enriques surfaces and other
non-Pfaffian subcanonical subschemes of
codimension 3}, Special issue in honor of Robin
Hartshorne, Comm. Algebra  28  (2000), 
5629-5653.

\item{[5]} D. Eisenbud, S. Popescu, C.
Walter, {\it Lagrangian subbundles and
codimension 3 subcanonical subschemes}, Duke
Math. J. 107 (2001), 427-467. 

\item{[6]} H. Grauert, G. M\"ulich, {\it
Vektorb\"undel vom Rang 2 \"uber dem
$n$-dimensionalen komplex-projektiven Raum},
Manusc. Math., 16 (1975), 75-100.

\item{[7]} R. Hartshorne, {\it Varieties of
small codimension in projective space}, Bull.
AMS 80(6) (1974), 1017-1032. 

\item{[8]} R. Hartshorne, {\it Algebraic
Geometry}, Springer 1977.


\item{[9]} G. Horrocks, {\it A construction
for locally free sheaves}, Topology 7 (1968),
117-120.

\item{[10]} J.J. Rotman, {\it An introduction to
cohomological algebra}, Academic Press, 1979.

\item{[11]} J.P. Serre. {\it Sur les modules
projectifs}, S\'eminaire Dubreil-Pisot
(1960/61), Secr. Math. Paris, expos\'e 2
(1961).

\item{[12]} J. A. Vogelaar, {\it
Constructing vector bundles from codimension-two
subvarieties}, PhD thesis. Leiden 1978.

\item{[13]} C. Walter, {\it Pfaffian subschemes},
J. Algebraic Geom. 5 (1996), 671-704.

\bigskip

\centerline{Departamento de \'Algebra}
\centerline{Facultad de de Ciencias Matem\'aticas}
\centerline{Universidad Complutense de Madrid}
\centerline{28040 Madrid, Spain}
\centerline{\tt arrondo@mat.ucm.es}

\end